\documentclass[12pt]{article}
\usepackage[T2A]{fontenc}
\usepackage[english]{babel}
\usepackage{amsfonts,amssymb,amsmath,enumerate,amsthm}
\usepackage{epsfig}
\usepackage{float}

\newcommand{\RR}{\mathbb R}

\newcommand{\NN}{{\mathbb N}}

\newtheorem{theorem}{Theorem}

\newtheorem{remark}{Remark}

\newtheorem{corollary}{Corollary}
\newtheorem{example}{Example}
\newtheorem{proposition}{Proposition}

\newcommand{\beq}{\begin{equation}}
\newcommand{\eeq}{\end{equation}}
\newcommand{\ba}{\begin{array}}
\newcommand{\ea}{\end{array}}
\newcommand{\bea}{\begin{eqnarray}}
\newcommand{\eea}{\end{eqnarray}}

\providecommand{\keywords}[1]
{
  \small	
  \textbf{Keywords: } #1
}
\providecommand{\AMS}[1]
{
  \small	
  \textbf{AMS Subject classification: } #1
}

\usepackage{graphicx}
\usepackage[usenames,dvipsnames]{color}
\usepackage{transparent}

\DeclareMathAlphabet{\mathpzc}{OT1}{pzc}{m}{it}

\usepackage{ulem}
\usepackage{cancel}

\begin{document}
\begin{center}

{\bf On reconstruction from imaginary part
for radiation solutions in two dimensions}
\vskip 10pt
{\it A.V. Nair,  \it R.G. Novikov }

\vskip 10pt

{\centering\footnotesize  \par}

\end{center}

\begin{abstract}
We consider a radiation solution $\psi$ for the Helmholtz equation in an exterior region in $\RR^2$.
We show that $\psi$ in the exterior region is uniquely determined by its imaginary part $Im(\psi)$ on an interval of a line $L$ lying in the exterior region. This result has a holographic prototype in the recent work Nair, Novikov (2025, J. Geom. Anal. 35, 123).
Some other curves for measurements instead of the lines $L$ are also considered. Applications to the Gelfand-Krein-Levitan inverse problem (from boundary values of the spectral measure in $\RR^2$) and to passive imaging are also indicated.
\end{abstract} 
\keywords{Two-dimensional Helmholtz equation, radiation solutions, Gelfand-Krein-Levitan inverse problem, passive imaging, holographic type uniqueness} \\
\AMS{35J05, 35J08, 35P25, 35R30}
\section{Introduction}
We consider the two-dimensional Helmholtz equation
\begin{equation}
\label{eq:1.1}
\Delta\psi(x) + \kappa^2\psi(x) = 0, \ \  x \in {\cal U},\ \  \kappa>0, 
\end{equation}
where $\Delta$ is the Laplacian in $x$, and $\cal U$ is an exterior region (open unbounded connected set) in $\RR^2$
consisting of all points outside a closed bounded regular curve $S$ (for example as in \cite{SK}). 
For equation (\ref{eq:1.1}) we consider the radiation solution $\psi$ such that:
$\psi$ is of class $C^2$ and satisfies the Sommerfeld's radiation condition
\begin{equation}
\label{eq:1.2}
\sqrt{|x|}\bigl(\frac{\partial}{\partial|x|} - i\kappa \bigr)\psi(x) \to 0 \ \ as \ \ |x| \to +\infty,
\end{equation}
uniformly in $x/|x|$.

Let
\begin{equation}
L = L_{x_0, \theta}=\{x\in\RR^2:\ \ x=x(s)=x_0+s\theta,\ \ -\infty<s <+\infty  \}, \nonumber
\end{equation}
\begin{equation}  
 L^+=L^+_{x_1, \theta}=\{x\in\RR^2:x=x(s)=x_1+s\theta,\ \ 0<s <+\infty  \}, \nonumber
\end{equation}
\begin{equation}
\label{eq:1.3}
 L^-=L^-_{x_2, \theta}=\{x\in\RR^2:\ \ x=x(s)=x_2-s\theta,\ \ 0<s <+\infty  \},
\end{equation}
where $x_0,x_1,x_2 \in \RR^2$, $\theta \in \mathbb S^1$ and $x_1,x_2 \in L$. \\ 
The present work combines results in \cite{NN}, and methods of \cite{NV} (the three-dimensional case), to obtain analogical results to \cite{NV} in the two-dimensional case. \\
In particular, we show that for any straight line $L\subset  {\cal U}$, any complex-valued radiation solution  $\psi$  on ${\cal U}$ is uniquely determined  by  the  $Im(\psi)$
on an arbitrary non-empty interval $\Lambda$ of $L$; see Theorem 2 (which is a corollary of Theorem 1) and Corollary 1 in Section 2. We also consider other curves for measurements instead of the lines $L$;
see Example 1, Theorem 3, and Remark 1 in Section 2. 

The aforementioned Theorems 1, 2, and Corollary 1 have holographic prototypes in \cite{NN}, where reconstruction of $\psi$ is considered from the intensity $|\psi_0 + \psi|^2$ in place of $Im (\psi)$. Here, $\psi_0=e^{ikx}$ is a plane wave solution of (\ref{eq:1.1}), i.e., $k\in \RR^2, |k| = \kappa$.
In particular, results of \cite{NN} solve one of the old mathematical questions of holography in its two-dimensional setting and admit straightforward applications to phaseless
inverse scattering in two dimensions. In turn, \cite{NN} continues studies of \cite{N}
(the three-dimensional
case).

In the present work, as in \cite{NV}, we are motivated by the Gelfand-Krein-
Levitan inverse problem (from boundary values of the spectral measure in
the whole space) and by passive imaging.
As relevant literature, see, e.g., \cite{AA}, \cite{AH}, \cite{YB}, \cite{BS}, \cite{GG}, \cite{GB}, \cite{BM}, \cite{M}, \cite{NV}, \cite{WL}.   In addition, the present work is devoted to the two-dimensional case.

The Gelfand-Krein-Levitan problem in question (in its fixed energy version in dimension $d = 2$) consists in determining the potential $v$ in the Schr{\"o}dinger equation 
\begin{equation}
\label{eq:1.4}
-\Delta\psi(x) + v(x)\psi(x) = \kappa^2\psi(x) + \delta(x-y), \ \ x, y \in \RR^2, \ \ \kappa > 0,
\end{equation}
from the imaginary part of its radiation solutions $\psi=R_v^+(x,y, \kappa)$ (that is, satisfying (\ref{eq:1.2})) for one $\kappa$ and all $x, y$ on some part of the boundary of a domain containing the support of $v$. Here, $\delta$ is the Dirac delta function. In this problem, Im$(R_v^+)$ is related to the spectral measure (in spectral decomposition) of the Schr{\"o}dinger operator $H = -\Delta + v$ in $L^2(\RR^2)$, at least, for real-valued $v$; see, \cite{YB}.
For introduction into the spectral theory of the Schrödinger operator, see, e.g., \cite{BSH}, \cite{H}. 

In applications, the problem of recovering $v$ from
boundary values of $Im (R^+_{v})$  at fixed $\kappa$ arise, in particular, in the framework of passive acoustic tomography (in ultrasonics, ocean acoustics, local helioseismology).
In these frameworks, equation (\ref{eq:1.4}) is considered as the acoustic Helmholtz equation, where complex-valued $v(x) = v(x, \kappa)$ is related to the perturbation of the refraction index
and $Im (R^+_{v})$ is related to cross correlations of wave fields generated by random sources; see, e.g., formula (49) in \cite{GB}.

For more information about the aforementioned Gelfand-Krein-Levitan problem and its relevance to passive imaging, see \cite{AA}, \cite{YB}, \cite{NV} and references therein.

Theorems 1,  2, \ Corollary 1, and Theorem 3 mentioned above admit the same applications to the Gelfand-Krein-Levitan inverse problem
and passive imaging  in two dimensions as their prototypes in \cite{NV} in three dimensions.
In particular, in the present work we give global uniqueness results for the aforementioned  Gelfand-Krein-Levitan inverse problem
at fixed energy in dimension $d=2$; see Theorems 4,  5, and Remark 1 in Subsection 2.2.
These results also contribute to two-dimensional inverse scattering at fixed energy, in general; as relevant literature on this topic,
see, e.g., \cite{AB},  \cite{PGG}, \cite{NR}, \cite{GN},  \cite{RGN}, \cite{IAT}, and references therein.

In the present work, we use the Karp expansion (\ref{eq:3.2}) below for the radiation solutions $\psi$ of equation (\ref{eq:1.1}) instead of the Atkinson-Wilcox expansion, used in \cite{NV}, for radiation solutions of the Helmholtz equation in three dimensions. In addition, we use very recent results on the Karp expansion obtained in \cite{NN}.

% Suppose that
% \begin{equation}
% \label{eq:1.5}
% supp \ v \subset \RR^2\setminus\overline{\cal{U}}, \ \  y \in \RR^2\setminus\overline{\cal{U}},\ \ \rho > 0,
% \end{equation}
% where $\overline{\cal{U}}$ is the closure of $\cal{U}$ (and $\cal{U}$ as in (\ref{eq:1.1})).
The main results of the present work are presented in more detail and proved in Section 2, Subsection 3.2, and Section 4.
In our proofs, we proceed from the results recalled in Section~3.

\section{Main results}
\subsection{Determination of radiation solutions $\psi$ from Im$(\psi)$}
Our key result is as follows.

\begin{theorem}\label{thm:1}
Let $\psi$  be a radiation solution of equation (\ref{eq:1.1}) as in  (\ref{eq:1.3}).
Let $L$, $L^+$ and $L^-$ be as given in (\ref{eq:1.3}) such that  $L^+=L^+_{x_1, \theta}\subset {\cal U}$ and $L^-=L^-_{x_2, \theta}\subset {\cal U}$, $x_1,x_2 \in L$, where ${\cal U}$ is the region in  (\ref{eq:1.1}). Then $\psi$ on $L^+$ $\cup$ $L^-$ is  uniquely determined by Im$(\psi)$ on   $\Lambda^+\cup \Lambda^-$,
where $\Lambda^+$ and $\Lambda^-$ are arbitrary non-empty intervals of $L^+$ and $L^-$, respectively.
\end{theorem}

As a corollary, we also get the following result.

\begin{theorem}\label{thm:2}
Let $\psi$  be a radiation solution of equation (\ref{eq:1.1}) as in  (\ref{eq:1.3}).
Let $L$ be as given in (\ref{eq:1.3}) such that  $L=L_{x_0, \theta}\subset {\cal U}$, where ${\cal U}$ is the region in  (\ref{eq:1.1}). Then $\psi$ on $L$ is  uniquely determined by Im$(\psi)$ on   $\Lambda$,
where $\Lambda$ is an arbitrary non-empty interval of $L$.
\end{theorem}
Theorem 1 is proved in Section 4 using the Karp expansion of \cite{SK} for the radiation solutions of equation (\ref{eq:1.1}), results of \cite{NN}, and methods of \cite{NV}. Similar to \cite{NV}, we use a two-point approximation for $\psi$ in terms of $Im(\psi)$; see Proposition 1 in Section 3.

Note that $\psi$ and $Im(\psi)$ are real-analytic on $\cal{U}$, and therefore on $L^+ \cup L^-$ in Theorem 1 \ and on $L$ in Theorem 2. Because of this analyticity, Theorem 1 reduces to the case when $\Lambda^+\cup \Lambda^- = L^+ \cup L^-$ and Theorem 2 reduces to the case when $\Lambda = L$.

Theorem 2 is proved as follows (for example). We assume that $\Lambda=L$. Then we simply consider $L^+ \subset L$ and $L^-\subset L$ such that $L^+ \cap L^- \neq \emptyset$. Then we may apply Theorem 1. 

\begin{corollary}
Under the assumptions of Theorem 2, Im$(\psi)$ on  $L$ uniquely determines $\psi$  
in the entire region ${\cal U}$.  

\end{corollary}
Corollary 1 follows from  Theorem 2, formula (\ref{eq:3.9}) recalled in Section 3, and analyticity of  $\psi$ in ${\cal U}$.

Theorems 1 and 2, and Corollary 1 have holographic prototypes in \cite{NN}. \\

In connection with other curves of measurements instead of the lines $L$ our results are as follows. 

The results in Theorems 1 and 2, don't hold for some other curves in place of $L$. An example is as follows.
Let
\begin{equation}
\label{eq:2.1}
\mathbb{S}_r = \{x \in \RR^2: |x| = r\}, \ \  r > 0.
\end{equation}
\begin{example}
Let $\psi(x) = G^+(x,\kappa) = \frac{i}{4}H_0(\kappa|x|)$, where $H_0$ is the Hankel function of first kind of order zero. Then $\psi$ is a non-zero radiation solution of equation (\ref{eq:1.1}) if $\{0\} \subset \RR^2 \setminus \overline{\cal{U}}$, but $Im(\psi) \equiv 0$ on the circles $\mathbb{S}_r$ for $r = \frac{c_j}{\kappa}, \ j \in \NN$, where $c_j$ are the positive real roots of the Bessel function $J_0 = Re(H_0)$. 
\end{example}
The fact that there are infinite positive real roots of $J_0$ and they tend to positive infinity can be found in \cite{AS}.
Nevertheless, as in the three-dimensional scenario, the uniqueness results of Theorem 2 and Corollary 1, remain valid for many interesting curves instead of the lines $L$.

Suppose that
\begin{equation}
\label{eq:2.2}
\begin{split}
\cal{D} \text{ is an open bounded domain in  $\RR^2$},\\
\Gamma = \partial \cal{D} \text{ is real analytic and connected.}
\end{split}
\end{equation}

\begin{theorem}
Let $\psi$ be a radiation solution of equation (\ref{eq:1.1}) as in (\ref{eq:1.2}). Let $\Gamma$ be a curve as above, where $\kappa$ is not a Dirichlet eigenvalue for $\cal{D}$, and $\overline{\cal{D}} = \cal{D} \ \cup$ $\Gamma$ $\subset \cal{U}$. Then $\psi$ on $\cal{U}$ is uniquely determined by Im$(\psi)$ on $\Lambda$ where $\Lambda$ is an arbitrary non-empty open interval of $\Gamma$.   
\end{theorem}
The proof of Theorem 3 is similar to that of Theorem 4 in \cite{NV}. In particular, this proof uses real analytic continuations, solving the Dirichlet problem for the Helmholtz equation in $\cal{D}$, and Corollary 1 in place of its three-dimensional prototype in \cite{NV}.
\subsection{Applications to Gelfand-Krein-Levitan problem}

The aforementioned results on recovering a radiation solution $\psi$ from Im$(\psi)$ give a reduction of the Gelfand-Krein-Levitan problem (of inverse spectral theory and passive imaging in dimension $d = 2$) to the inverse scattering problem of finding $v$ in (\ref{eq:1.4}) from boundary values of $\psi = R_v^+$.
This reduction and known results imply, in particular, that $v$ in (\ref{eq:1.4}) is uniquely determined by Im$(R_v^+(x, y, \kappa))$ for one $\kappa$ and all $x, y$ on an arbitrary open interval $\Lambda$ of $L$, where $supp \ v \subset \Omega$, $L$ is a line in $\RR^2\setminus\overline{\Omega}$, $\Omega$ is an open bounded connected domain in $\RR^2$, $\overline{\Omega}$ is the closure of $\Omega$; see Theorem 4 below.

\begin{theorem}
Let $v \in L^\infty(\RR^2)$, $supp \  v \subset \Omega$, and $L \subset \RR^2\setminus\overline{\Omega}$, where $\Omega$ is an open bounded connected domain in $\RR^2$, $\overline{\Omega}$ is the closure of $\Omega$, and $L$ is a straight line. Let $R_v^+(x, y, \kappa)$ be the outgoing Green function for equation (\ref{eq:1.4}). Then $v$ is uniquely determined by Im$(R_v^+(x, y, \kappa))$ for one $\kappa$ and all $x, y$ on $\Lambda$, where $\Lambda$ is an arbitrary non-empty open interval of $L$.
\end{theorem}

\begin{theorem}
Let $v\in L^\infty(\RR^2), \  supp \ v \subset \RR^2\setminus\overline{\cal{U}}$, where $\cal{U}$ is as in (\ref{eq:1.1}), and property (\ref{eq:2.4}) hold for fixed $\kappa > 0$.
Let $\Gamma$ be a curve as in (\ref{eq:2.2}), where $\kappa$ is not a Dirichlet eigenvalue for $\cal{D}$, and $\overline{\cal{D}} = \cal{D} \ \cup$ $\Gamma$ $\subset \cal{U}$. Then $v$ is uniquely determined by Im$(R_v^+(x,y,\kappa))$ on $\Lambda \times \Lambda$, where $\Lambda$ is an arbitrary non-empty open interval of $\Gamma$.
\end{theorem}
In Theorems 4 and 5, we do not assume that $v$ is real-valued, but we assume that (\ref{eq:2.4}) formulated in Subsection 3.4 holds.

Theorems 4 and 5 can be proved in a similar way to Theorems 3 and 5 in \cite{NV}, respectively.  These proofs use Theorems 1,2, and 3 in place of their three-dimensional prototypes in
\cite{NV}, some known results on two-dimensional direct scattering in place of similar results on three-dimensional direct scattering used in \cite{NV},
and also the two-dimensional global uniqueness result on recovering from the scattering amplitude $A$  in (\ref{eq:3.14}) at fixed $\kappa$ (see  Corollary 5.1 in \cite{AB} in addition to results of \cite{NR}, \cite{GN}).
These proofs follow from these indications in a 
straightforward 
way.

\begin{remark}
    The case is also of interest when the boundary $\Gamma$ in (\ref{eq:2.2}) is not connected but consists of two disjoint connected components $\Gamma_1$ and $\Gamma_2$, where $\Gamma_1=\partial {\cal D}_1$, $\Gamma_2=\partial {\cal D}_2$ and  ${\cal D}_1$,  ${\cal D}_2$ are open bounded domains such that $\RR^2 \setminus {\cal U} \subset {\cal D}_1  \subset {\cal D}_2$.
In this case Theorem 3 is valid with $\Lambda$ replaced by $\Lambda_1 \cup \Lambda_2$, whereas Theorem 4 is valid with $\Lambda \times \Lambda$ replaced by $(\Lambda_1 \cup \Lambda_2) \times \Lambda_1$
(as well as with $\Lambda \times \Lambda$ replaced by $(\Lambda_1 \cup \Lambda_2) \times \Lambda_2$, where $\Lambda_1$, $\Lambda_2$ are arbitrary non-empty open intervals of $\Gamma_1$ and $\Gamma_2$, respectively.
\end{remark}
Remark 1 is similar to Remark 8 in \cite{NV} and is motivated by Theorem 1 in \cite{AH} and Theorem 4.3 in \cite{BM}, where uniqueness results for monochromatic passive imaging in three dimensions are given for the case of a boundary with two disjoint connected components.

\section{Preliminaries}

Let $(r,\phi)$ denote the polar coordinates of a point $x$ and let $\theta = \frac{x}{|x|}$ be its direction. So $-\theta$ will be the direction of $-x$ whose polar coordinates will be $(r, \phi + \pi)$. \\
Let
\beq
\label{eq:3.1}
B_{\rho}=\{x\in\RR^2 \ | \ |x|<\rho\},\ \ \rho>0. 
\eeq
We use the following asymptotic expansion for radiation solution $\psi$ of equation (\ref{eq:1.1}):
\begin{equation}
\label{eq:3.3}  
\psi(x)\sim \sqrt{\frac{2}{\pi \kappa |x|}}e^{i(\kappa |x| - \frac{\pi}{4})}\sum\limits_{j=0}^{\infty}\frac{f_j(\phi)}{|x|^{j}} \ \ \mbox{\rm for}\ \ x\in \RR^2, \ \ |x| \to \infty.
\end{equation}
Expansion (\ref{eq:3.3}) looks similar to the convergent Atkinson-Wilcox expansion in three dimensions (in \cite{A}).
However, the series in (\ref{eq:3.3}) diverges in general and in particular for $\psi(x) = H_0(\kappa |x|)$ (see \cite{SK} for more details).  \\
\subsection{The Karp Expansion}

Due to \cite{SK},  we have the following Karp expansion: 

Suppose that  $\psi$ is a radiation solution of  equation (\ref{eq:1.1}), and  $\RR^2\setminus B_{\rho}  \subset {\cal U}$. \\

\begin{equation}
\label{eq:3.2}  
\psi(x)=H_0(\kappa |x|)\sum\limits_{j=0}^{\infty}\frac{F_j(\phi)}{|x|^{j}} + H_1(\kappa |x|)\sum\limits_{j=0}^{\infty}\frac{G_j(\phi)}{|x|^{j}},
\end{equation}
\begin{center}
{$\mbox{\rm for}\ \ x\in \RR^2\setminus B_\rho,$}
\end{center}
where $H_0$ and $H_1$ are the Hankel functions of the first kind of order zero and one respectively. The series converges absolutely and uniformly for $|x| \geqq \rho_1 > \rho$ (for any such $\rho_1$). 

The important result is that $F_j(\phi), G_j(\phi), \ j = 0,1,2,...n,$ in (\ref{eq:3.2}) are uniquely determined from $f_j(\phi), f_j(\phi + \pi), \ j = 0,1,2,...n,$ in (\ref{eq:3.3}) via recurrent formulas; (see Theorem III of \cite{SK} for $n=0$, Theorem 3 of \cite{NN} for the general case).

\subsection{A two-point approximation for $\psi$}
Let
\begin{equation}
\label{eq:3.4}
I(x) = \sqrt{|x|}Im(\psi(x)), \ \ x \in \cal{U},
\end{equation}
where $\psi$ is a radiation solution of equation (\ref{eq:1.1}). \\

We have that
\begin{equation}
\label{eq:3.5}
2iI(x) = \sqrt{\frac{2}{\pi \kappa}}(e^{i(\kappa |x| - \frac{\pi}{4})}f_0(\phi) - e^{-i(\kappa |x| - \frac{\pi}{4})}\overline{f_0(\phi)}) + O\bigl(\frac{1}{|x|}\bigr), \text{ as } |x| \to \infty,
\end{equation}
uniformly in $\theta$ = $x/|x|$, where $f_0$ is the leading coefficient in (\ref{eq:3.3}).\\
\begin{proposition}
Let $\psi$ be a radiation solution of equation (\ref{eq:1.1}). Then
\begin{equation}
\label{eq:3.6}
f_0(\phi) = \frac{1}{sin(\kappa \tau)}\sqrt{\frac{\pi \kappa}{2}}\bigl(I(x)e^{-i(\kappa |y| - \frac{\pi}{4})} - I(y)e^{-i(\kappa |x| - \frac{\pi}{4})} + O\bigl(\frac{1}{|x|}\bigr) \bigr),  \text{ as } |x| \to \infty,
\end{equation}
\begin{center}
$x,y \in L_{x_0,\theta}^+$,  $x_0 = 0$,  $y = x + \tau\theta$,  $\theta \in \mathbb{S}$,  $\tau > 0$,
\end{center}
uniformly in $\theta$, where $f_0$ is the leading coefficient in (\ref{eq:3.3}), I is defined in (\ref{eq:3.4}), $L^+$ is the ray defined by (\ref{eq:1.3}), and $sin(\kappa\tau)\neq0$ for fixed $\tau$.
\end{proposition}
Formula (\ref{eq:3.6}) is a two-point approximation for $f_0$ and together with (\ref{eq:3.3}) also gives a two-point approximation for $\psi$ in terms of $I$. In phaseless inverse scattering and holography, formulas of such a type go back to \cite{N1} (see also \cite{NS}, \cite{NSh}, \cite{N}, \cite{NN}, \cite{NV}).
\begin{remark}
If an arbitrary function $I$ on $L_{0,\theta}$ satisfies (\ref{eq:3.5}), then formula (\ref{eq:3.6}) holds, for fixed $\theta \in \mathbb{S}$.
\end{remark}
We obtain (\ref{eq:3.6}) from the system of equations for $f_0, \overline{f_0}$:
\begin{equation}
\label{eq:3.7}
\sqrt{\frac{2}{\pi \kappa}}(e^{i(\kappa |x| - \frac{\pi}{4})}f_0(\phi) - e^{-i(\kappa |x| - \frac{\pi}{4})}\overline{f_0(\phi)}) = 2iI(x) + O(|x|^{-1}),
\end{equation}
\begin{equation}
\sqrt{\frac{2}{\pi \kappa}}(e^{i(\kappa |y| - \frac{\pi}{4})}f_0(\phi) - e^{-i(\kappa |y| - \frac{\pi}{4})}\overline{f_0(\phi)}) = 2iI(y) + O(|y|^{-1}),  \nonumber
\end{equation}
where $x,y$ are as in (\ref{eq:3.6}). In particular, we use that $|y| = |x| + \tau$. 
In turn, (\ref{eq:3.7}) follows from (\ref{eq:3.5}).
\subsection{A Green type formula}
The following formula holds:
\begin{align}
&\psi(x)=  2\int_{L}\frac{\partial G^+(x-y,\kappa)}{\partial \nu_{y}}\psi(y)dy,\ \ x\in V_L,  \label{eq:3.9}\\
&G^+(x,\kappa)= \frac{i}{4}H_0(\kappa|x|),\ \ x\in \RR^2, \nonumber
\end{align}
where $\psi$ is a radiation solution of equation (\ref{eq:1.1}), $L$ and $V_L$ are line and open half-plane in ${\cal U}$,
where $L$ is the boundary of $V_L$, $\nu$  is the outward normal to $L$ relative to $V_L$; see, for example, formula 5.84 in \cite{BR}.

\subsection{Some facts of direct scattering}
The radiation solution $R_v^+ = R_v^+(.,y, \kappa)$ for equation (\ref{eq:1.4}) satisfies the integral equation 
\begin{equation}
\label{eq:2.3}
R_v^+(x,y, \kappa) = -G^+(x-y, \kappa) + \int_{\Omega}G^+(x-y, \kappa)v(z)R_v^+(z,y, \kappa)dz,
\end{equation}
where $x,y \in \RR^2$, $G^+$ is given in (\ref{eq:3.9}). \\
We consider equations (\ref{eq:1.4}) and (\ref{eq:2.3}) assuming for simplicity that

\begin{equation}
\label{eq:3.10}
    \text{$v$ is complex-valued}, \ v\in L^\infty(\Omega), \  v \equiv 0 \text{ on } \RR^2\setminus\overline{\Omega}, \  
\end{equation}
and
\begin{equation}
\label{eq:2.4}
\text{equation (\ref{eq:2.3}) is uniquely solvable for } R_v^+(.,y, \kappa) \in L^2(\Omega). 
\end{equation}

It is known that if $v$ satisfies (\ref{eq:3.10}) and is real-valued (or Im$(v) \leq 0$) then (\ref{eq:2.4}) is fulfilled automatically; see, for example, \cite{CK}, at least in three dimensions.

Note that
\begin{equation}
\label{eq:3.12}
    R_v^+(x,y, \kappa) = R_v^+(y,x, \kappa), \ \ x, y \in \RR^2.
\end{equation}

We also consider the scattering wave functions $\psi^+$ for the homogeneous equation (\ref{eq:1.4}) (i.e without $\delta$):
\begin{equation}
\label{eq:3.11}
    \psi^+ = \psi^+(x. k) = e^{ikx} + \psi^+_{sc}(x, k), \ \ x,k \in \RR^2, \ |k| = \kappa,
\end{equation}
where $\psi^+_{sc}(x, k)$ satisfies the radiation condition (\ref{eq:1.2}) at fixed $k$. 

The following formulas hold:
\begin{equation}
\label{eq:3.13}
    R_v^+(x,y, \kappa) = -\frac{1}{2}\sqrt{\frac{1}{2\pi\kappa|x|}}e^{i(\kappa|x|+\frac{\pi}{4})}\psi^+(y, -\kappa\frac{x}{|x|}) + O\bigl(\frac{1}{|x|^{3/2}}\bigr), \ \text{ as } |x| \to \infty,
\end{equation}
at fixed $y$;
\begin{equation}
\label{eq:3.14}
   \psi^+_{sc}(x, k) = \frac{e^{i\kappa|x|}}{|x|^{1/2}}A(k, \kappa \frac{x}{|x|}) + O\bigl(\frac{1}{|x|^{3/2}}\bigr), \ \text{ as } |x| \to \infty,
\end{equation}
at fixed $k$, where $A$ arising in (\ref{eq:3.14}) is the scattering amplitude for the homogeneous equation (\ref{eq:1.4}).\\
In connection with the aforementioned facts concerning $R^+_v$ and $\psi^+$, see, e.g., \cite{FM}.
\section{Proof of Theorem 1}
\subsection{Case $L^+_{x_1,\theta} \cup L^-_{x_2,\theta} \subseteq L_{0,\theta}$}
In a similar way to \cite{NV} (where the three-dimensional case is considered), we prove that Im$(\psi)$ on $L^+_{x_1,\theta} \cup L^-_{x_2,\theta} \subseteq L_{0,\theta}$ defined as in (\ref{eq:1.3}), where ${0}$ denotes the origin in $\RR^2$, uniquely determines  $f_j(\phi)$  in (\ref{eq:3.3}) $\forall j \in \mathbb{N}$. The proof is given below. Then, using formulas in \cite{SK} and \cite{NN} as recalled in Subsection 3.1, all $F_j(\phi), G_j(\phi),$ in (\ref{eq:3.2}) are uniquely determined by $f_j(\phi), f_j(\phi + \pi)$. The rest follows from the convergence of the series in (\ref{eq:3.2}) and analyticity of $\psi$  and Im$(\psi)$ on $L^+\subseteq L^+_{0, \theta}$.

The determination of $f_0(\phi)$ from Im$(\psi)$ on $L^+_{x_1,\theta}$ follows from (\ref{eq:3.6}). 

Suppose that $f_0(\phi),f_1(\phi),f_2(\phi),...f_n(\phi)$ are determined then the determination of $f_{n+1}(\phi)$ from Im$(\psi)$ on $L^+_{x_1,\theta}$ is as follows.

Let
\begin{equation}
\label{eq:4.4}  
\psi_{n}(x)=\sqrt{\frac{2}{\pi \kappa |x|}}e^{i(\kappa |x| - \frac{\pi}{4})}\sum\limits_{j=0}^{n}\frac{f_j(\phi)}{|x|^{j}}, \ \ \mbox{\rm where}\ \ \theta=\frac{x}{|x|} = (cos(\phi), sin(\phi)),
\end{equation}
\begin{equation}
\label{eq:4.5}  
I_n(x)=\sqrt{|x|}Im(\psi(x)),
\end{equation}
\begin{equation}
\label{eq:4.6}  
J_n(x)=|x|^{n+1}(I(x)-I_n(x)),
\end{equation}
where $x$ is as in (\ref{eq:3.2}), $I(x)$ is defined by (\ref{eq:3.4}).

We have that:
\begin{align}
&2iI(x) = 2I_n(x) + \sqrt{\frac{2}{\pi \kappa}}\frac{e^{i(\kappa |x| - \frac{\pi}{4})}}{|x|^{n+1}}f_{n+1}(\phi) - \sqrt{\frac{2}{\pi \kappa}}\frac{e^{-i(\kappa |x| - \frac{\pi}{4})}}{|x|^{n+1}}\overline{f_{n+1}(\phi)}  \label{eq:4.7}\\
&+O\bigr(\frac{1}{|x|^{n+2}}\bigl), \nonumber
\end{align}
\begin{equation}
\label{eq:4.8}  
2iJ_n(x) = \sqrt{\frac{2}{\pi \kappa}}(e^{i(\kappa |x| - \frac{\pi}{4})}f_{n+1}(\phi) - e^{-i(\kappa |x| - \frac{\pi}{4})}\overline{f_{n+1}(\phi)}) + O\bigl(\frac{1}{|x|}\bigr),
\end{equation}
as $|x|\to +\infty$, uniformly in $\theta=x/|x|$, where $I$ is defined by (\ref{eq:3.4}).

Due to (\ref{eq:4.8}) and Remark 2, as $|x| \to +\infty$ we get:
\begin{align}
&f_{n+1}(\phi) = \frac{1}{sin(\kappa \tau)}\sqrt{\frac{\pi \kappa}{2}}\bigl(J_n(x)e^{-i(\kappa |y| - \frac{\pi}{4})} - J_n(y)e^{-i(\kappa |x| - \frac{\pi}{4})} + O\bigl(\frac{1}{|x|}\bigr)\bigr) ,\label{eq:4.9}\\
&x,y\in L^+\subseteq L^+_{0, \theta},\ \  y=x+\tau\theta,\ \  \theta\in\mathbb S^1,\ \  \tau>0,  \nonumber
\end{align}
uniformly in $\theta$, assuming that  $sin(\kappa \tau) \neq 0$ for fixed $\tau$ (where the parameter $\tau$ can always be fixed in such a way for fixed $\kappa > 0$).

Formulas (\ref{eq:3.4}), (\ref{eq:4.4})-(\ref{eq:4.6})  and (\ref{eq:4.9}) determine  $f_{n+1}(\phi)$, give the step of induction for finding all $f_j(\phi)$ from Im$(\psi)$ on $L^+_{x_1, \theta}$. The determination of $f_j(\phi + \pi)$ from Im$(\psi)$ on $L^-_{x_2,\theta}$ is completely similar.

This completes the proof of Theorem 1 for the case $L^+_{x_1,\theta} \cup L^-_{x_2,\theta} \subseteq L_{0,\theta}$.

\subsection{General Case}
In fact, in a similar way with \cite{N}, \cite{NV}, and \cite{NN}, the general case reduces to the case of Subsection 4.1 by the change of variables 
\begin{align}
&x'=x-q \label{eq:6.1}\\
&\mbox{\rm for some fixed }\ \ q\in\RR^2\ \ \mbox{\rm such that}\ \ L^+\subseteq L^+_{q, \theta},\ \  L^-\subseteq L^-_{q, \theta}.   \nonumber
\end{align}
Here, without restriction of generality we assume that $L^+ \cap L^- = \emptyset$.
In the new variables $x'$, where $x'\in{\cal U'}={\cal U}-q$ we have that:
\begin{equation}
\label{eq:6.5}  
L^+\subseteq L^+_{q, \theta}=L^+_{0, \theta} \;\;\;\;\;\;\;
L^-\subseteq L^-_{q, \theta}=L^-_{0, \theta}.   
\end{equation}
In addition,
\begin{equation}
\label{eq:6.6}  
\psi(x')=\sqrt{\frac{2}{\pi \kappa |x'|}}e^{i(\kappa |x'| - \frac{\pi}{4})}\sum\limits_{j=0}^{\infty}\frac{f_j'(\phi)}{|x'|^{j}}  \ \ \mbox{\rm for}\ \ x'\in \RR^2\setminus B_{\rho'}, \ \  
\end{equation}
for some new $f_j'$, where $\rho'$ is such that $\RR^2\setminus B_{\rho'} \subset \cal{U'}$. \\
We also have
\begin{equation} 
\label{eq:6.7} 
\psi(x')=H_0(\kappa |x'|)\sum\limits_{j=0}^{\infty}\frac{F_j'(\phi)}{|x'|^{j}} + H_1(\kappa |x'|)\sum\limits_{j=0}^{\infty}\frac{G_j'(\phi)}{|x'|^{j}},
\end{equation}
where $F_j'$ and $G_j'$ are appropriate new coefficients (in accordance with subsection 6.2 of \cite{NN}), and the series converges absolutely and uniformly in $|x'| \geq r > \rho'$.

In view of (\ref{eq:6.6}), (\ref{eq:6.7}), we complete the proof of Theorem 1 by repeating the proof of Subsection 4.1.

\begin{remark}
Our proof of Theorem 1 has holographic prototypes in \cite{NN}. Additional formulas for finding $f_j$ can also be obtained using approaches of \cite{N4}, \cite{NS}, \cite{NG}.
\end{remark}

\pagebreak

\noindent
Arjun V. Nair\\
School of Mathematics,\\
IISER Thiruvanthapuram, Thiruvananthapuram, India;\\
E-mail: arjunnair172920@iisertvm.ac.in
\\
\\
\\

\noindent
Roman G. Novikov\\
CMAP, CNRS, Ecole Polytechnique,\\
Institut Polytechnique de Paris, 91128 Palaiseau, France;\\
IEPT RAS, 117997 Moscow, Russia;\\
E-mail: novikov@cmap.polytechnique.fr

\end{document}